\documentclass{amsproc}
\usepackage{amsmath, amssymb, amsthm,
amsfonts, amsthm, bm, eucal, graphicx}

\theoremstyle{plain}
\newtheorem{thm}{Theorem}[section]
\newtheorem{cor}[thm]{Corollary}
\newtheorem{lem}[thm]{Lemma}
\newtheorem{prop}[thm]{Proposition}

\theoremstyle{definition}
\newtheorem{defi}[thm]{Definition}
\newtheorem{conj}[thm]{Conjecture}
\newtheorem{nota}[thm]{Notation}
\newtheorem{obs}[thm]{Observation}
\newtheorem{obss}[thm]{Observations}
\newtheorem{rem}[thm]{Remark}
\newtheorem{rems}[thm]{Remarks}
\newtheorem{exa}[thm]{Example}
\newtheorem{exas}[thm]{Examples}
\newtheorem{prob}[thm]{Problem}
\newtheorem{sit}[thm]{}

\newcommand{\brem}{\begin{rem}}
\newcommand{\brems}{\begin{rems}}
\newcommand{\erem}{\end{rem}}
\newcommand{\erems}{\end{rems}}
\newcommand{\bexa}{\begin{exa}}
\newcommand{\bexas}{\begin{exas}}
\newcommand{\eexa}{\end{exa}}
\newcommand{\eexas}{\end{exas}}
\newcommand{\bdefi}{\begin{defi}}
\newcommand{\edefi}{\end{defi}}
\newcommand{\bdefis}{\begin{defis}}
\newcommand{\edefis}{\end{defis}}
\newcommand{\bcor}{\begin{cor}}
\newcommand{\ecor}{\end{cor}}
\newcommand{\blem}{\begin{lem}}
\newcommand{\elem}{\end{lem}}
\newcommand{\bconv}{\begin{conv}}
\newcommand{\econv}{\end{conv}}
\newcommand{\bconj}{\begin{conj}}
\newcommand{\econj}{\end{conj}}
\newcommand{\bprop}{\begin{prop}}
\newcommand{\eprop}{\end{prop}}
\newcommand{\bprob}{\begin{prob}}
\newcommand{\eprob}{\end{prob}}
\newcommand{\bthm}{\begin{thm}}
\newcommand{\ethm}{\end{thm}}
\newcommand{\bnota}{\begin{nota}}
\newcommand{\enota}{\end{nota}}
\newcommand{\bobs}{\begin{obs}}
\newcommand{\eobs}{\end{obs}}
\newcommand{\bobss}{\begin{obss}}
\newcommand{\eobss}{\end{obss}}
\newcommand{\bsit}{\begin{sit}}
\newcommand{\esit}{\end{sit}}
\newcommand{\be}{\begin{eqnarray}}
\newcommand{\ee}{\end{eqnarray}}
\newcommand{\bproof}{\begin{proof}}
\newcommand{\eproof}{\end{proof}}

\def\ba{\begin{array}}
\def\ea{\end{array}}
\def\bnum{\begin{enumerate}}
\def\enum{\end{enumerate}}

\newcommand{\ka}{{\mathbb K}}
\newcommand{\F}{{\mathbb F}}
\newcommand{\A}{{\mathbb A}}
\newcommand{\cO}{{\mathcal O}}

\newcommand{\N}{{\mathbb N}}

\newcommand{\T}{{\mathbb T}}

\newcommand{\Z}{{\mathbb Z}}
\newcommand{\cF}{{\mathcal F}}
\renewcommand{\phi}{\varphi}

\newcommand{\G}{{\Gamma}}

\newcommand{\ord}{{\operatorname{ord}}}

\newcommand{\diag}{{\operatorname{diag}}}

\newcommand{\ind}{{\operatorname{ind}}}
\newcommand{\End}{{\operatorname{End}}}
\newcommand{\rank}{{\operatorname{rank}}}
\newcommand{\Char}{{\operatorname{Char}}}

\newcommand{\Id}{{\operatorname{Id}}}
\newcommand{\Span}{{\operatorname{span}}}
\newcommand{\supp}{{\operatorname{supp}}}
\newcommand{\eval}{{\operatorname{eval}}}
\newcommand{\spec}{{\operatorname{spec}}}
\newcommand{\card}{{\operatorname{card}}}

\newcommand{\vecti}{{\operatorname{fragm}}}
\newcommand{\Vecti}{{\operatorname{Fragm}}}

\newcommand{\GF}{{\operatorname{GF}}}
\newcommand{\GL}{{\operatorname{GL}}}
\newcommand{\Tr}{{\operatorname{Tr}}}

\title[Discrete convolution operators]
{Discrete convolution operators in positive characteristic:\\
a variation on the Floquet-Bloch Theory}

\author{Mikhail Zaidenberg}
\address{Université
Grenoble I, Institut Fourier, UMR 5582 CNRS-UJF, BP 74, 38402 St.\
Martin d'H\`eres c\'edex, France} \email{zaidenbe@ujf-grenoble.fr}

\thanks{
\mbox{\hspace{11pt}}{\it 2010 Mathematics Subject
Classification}: 11C99, 11T99, 31C20, 37B15, 43A60.\\
\mbox{\hspace{11pt}}{\it Key words}: convolution operator,
lattice, finite field, discrete Fourier transform, discrete
harmonic function, pluri-periodic function.}

\thanks{{\bf Acknowledgements:} This is my pleasure to
thank Peter Kuchment for informing in his work and useful
discussions.}

\begin{document}

\begin{abstract}
The classical Floquet theory deals with Floquet-Bloch solutions of
periodic PDEs, see e.g. \cite{Ku3}. A discrete version of this
theory for difference vector equations on lattices, including the
Floquet theory on infinite periodic graphs, was developed by Peter
Kuchment \cite{Ku1, Ku2}. Here we propose a variation on this
theory for matrix convolution operators acting on vector functions
on lattices with values in a field of positive characteristic.
\end{abstract}

\maketitle
\date{}

\section{Introduction}
The classical Floquet theory founded in the work of Hill, Floquet,
Lyapunov, and Felix Bloch, deals with different types of periodic
ordinary differential equations. An advanced version of the theory
applies to the periodic PDEs and deals with the Floquet-Bloch
solutions \cite{Ku3}. A discrete analog appeals to the notion of a
periodic graph, see e.g. \cite{Co, vB}. This is an infinite graph
equipped with a free action of a free abelian group (a lattice)
$\Lambda$ of finite rank. Peter Kuchment \cite{Ku1}-\cite{Ku2}
adopted the Floquet theory for a periodic difference operator
$\Delta$ (e.g., the Laplace, Schr\"odinger, Markov operator, etc.)
acting in a translation invariant function space on a periodic
graph $\Gamma$; see also \cite{KuPi}, \cite{KuVa}. Assume that a
lattice $\Lambda$ acts on $\G$ with a finite number of orbits.
Then the following hold.
\begin{enumerate}\item Consider a   periodic
difference equation $\Delta(f)=0$ on an infinite periodic graph
$\G$. There exists a sequence $(f_n)$ of finite linear
combinations of Floquet solutions such that every other solution
$f$ can be decomposed in a series $f=\sum_{n=1}^\infty a_nf_n$
with constant coefficients.
\item Every $l_2$-solution can be approximated by solutions with
finite support.\item A criterion is provided as to when this
equation has no solution.\end{enumerate}

In these notes we propose a version of the discrete Floquet-Bloch
theory in positive characteristic. Our approach is purely
algebraic and ignores the analytic aspects, while preserving the
main features of the classical theory. In some respect, the
Floquet-Bloch theory for discrete operators gains a simplification
when passing to a positive characteristic. In particular, the
spectral theory occurs to be quite elementary. We consider matrix
convolution operators acting on vector functions on a lattice
$\Lambda$. Again, this includes the case of a periodic graph. As a
motivation we address the survey \cite{Za3}.

In positive characteristic, the traditional translation invariant
function spaces of the Floquet-Bloch theory (e.g., the span of all
Floquet functions) coincide with the space of all pluri-periodic
vector functions (see \ref{exas} and \ref{pre}). This issue, of
course, simplifies things. However, in positive characteristic the
latter space is not well adopted for applying the Fourier
transform\footnote{Instead of the Floquet transform as exploited
in \cite{Ku1}, \cite{Ku2}, we apply the usual Fourier transform.}.
More suitable in this respect is the subspace $\cF_p(\Lambda,
\ka)$, though it is characteristic-dependent. It consists of all
pluri-periodic vector functions on $\Lambda$ with values in the
base field $\ka$ whose period sublattices $\Lambda'\subseteq
\Lambda$ are of finite indices coprime to $p$. The role of Bloch
functions is played now by {\em elementary functions}. Such a
function is a product of a constant vector and a character of
$\Lambda$ with values in the multiplicative group $\ka^\times$ of
$\ka$.

Besides this Introduction, the paper contains sections 2-5. After
recalling in section 2 necessary preliminaries on harmonic
analysis on lattices in positive characteristic, we introduce our
principal function space and define the Fourier transform. The
rest of the paper is devoted to spectral analysis of convolution
operators. In section 3 we treat the scalar case, where such
operators are semi-simple (i.e. diagonalizable, see Propositions
\ref{scca} and \ref{sp1}). While in the vector case (see section
4) these operators are not any more semi-simple, in general. Given
a convolution operator $\Delta$ acting on vector functions, we
construct in Theorems \ref{kuch} and \ref{spde} a Jordan basis of
$\Delta$ consisting of elementary functions. Until subsection 4.3
the base field $\ka$ is assumed being algebraically closed. In
subsection 4.3  we deduce, for a generalized Laplace equation
$\Delta(f)=0$, the existence of a solution with values in a finite
field assuming that there exists such a solution with values in
its algebraic closure. In concluding section 5 we provide some
examples.

The author is grateful to the referee for correcting an error in a
formula in the first draft of the paper.

\section{Preliminaries}
\subsection{Harmonic analysis on lattices}
Let us recall some generalities; see e.g., \cite{Ni, Za2} for
details.

Given a prime integer $p>0$ we consider the finite Galois field
$\F_p=\GF(p)$ and its algebraic closure $\ka=\ka_p$. Recall that
the multiplicative group $\ka^\times$ of $\ka$ is a torsion group.
Moreover,  every finite subgroup of $\ka^\times$ is cyclic of
order coprime to $p$, and every cyclic group $\Z/n\Z$ of order
coprime to $p$ is isomorphic to a unique subgroup of $\ka^\times$;
see e.g., \cite{LN}.

\bsit\label{abr} {\em Lattices and dual tori.} Consider a lattice
$\Lambda$ of rank $s>0$. By a character of $\Lambda$ over $\ka$ we
mean a rational representation
$\chi:\Lambda\to\GL(1,\ka)\cong\ka^\times$. The image
$\chi(\Lambda)\subseteq \ka^\times$ is a cyclic subgroup of order
coprime to $p$. The  dual group $\Lambda^\vee$ of $\Lambda$ over
$\ka$ is the group of all characters $\Lambda\to\ka^\times$. It is
isomorphic to the algebraic torus $\T^s=(\ka^\times)^s$. By the
duality theorem $(\Lambda^\vee)^\vee\cong\Lambda$. Indeed, the
natural duality pairing $\Lambda\times\Lambda^\vee\to\ka^\times$
is non-degenerate.

More specifically, given a base $\mathcal V=(v_1,\ldots,v_s)$ of
$\Lambda$ we can identify $\Lambda$ with the standard integer
lattice $\Z^s$ by sending $v=\sum_{i=1}^s\lambda_iv_i\in\Lambda$
to $\lambda(v)=(\lambda_1,\ldots,\lambda_s)\in\Z^s$. The pairing
becomes  now
$$ \Z^s \times\T^s\to \ka^\times,\qquad (\lambda,z)\longmapsto
z^\lambda=\prod_{i=1}^s z_i^{\lambda_i}\quad\mbox{for}\quad
z=(z_1,\ldots,z_n)\in\T^s \,.$$ Fixing $z\in\T^s$ gives a
character $\Z^s\to\ka^\times$, $\lambda\longmapsto z^\lambda$, and
every character $\Z^s\to\ka^\times$ arises in this way. On the
other hand, every character of the torus $\T^s$ with values in
$\ka^\times$ is given by a Laurent monomial $z\longmapsto
z^\lambda$. This yields an identification
$(\Lambda^\vee)^\vee\cong\Lambda$. \esit

\bsit\label{coop} {\em Convolution operators.} Consider the group
algebra $\ka[\Lambda]$ over $\ka$ i.e., the algebra of functions
$\Lambda\to\ka$ with finite support endowed with the convolution
$$(f,g)\longmapsto f*g,\quad\mbox{where}\quad (f*g)(u)=\sum_{v\in \Lambda} f(v)
g(u-v)\,.$$ A convolution operator
$$\Delta_a \in\End(\ka[\Lambda]),\quad \Delta_a (f)
=a*f,\quad\mbox{where}\quad a\in\ka[\Lambda]\,,$$ can be written
as
$$\Delta_a=\sum_{v\in\Lambda} a(v)\tau_{-v}\,,$$
where $\tau_{v}$ denotes the shift by $v$:
 $\tau_{v}(f)(u)=f(u+v)$. Clearly, $\Delta_a$ commutes
with shifts: $[\Delta_a,\tau_{v}]=0$ $\forall v\in\Lambda$. In
fact, any endomorphism $\Delta\in\End_{\ka}\,(\ka[\Lambda])$
commuting with shifts is a convolution operator with kernel
$a=\Delta(\delta_0)$, where $\delta_v\in \ka[\Lambda]$ stands for
the delta-function concentrated on $v\in\Lambda$ i.e.,
$$\delta_v(v)=1\quad\mbox{and}\quad \delta_v(u)=0\,\,\,
\forall u\neq v\,.$$ For $a=\delta_v$ we obtain
$\Delta_a=\tau_{-v}$; indeed, $\delta_v *
\delta_0=\delta_v=\tau_{-v}(\delta_0)$. The $\ka$-algebra of
convolution operators is naturally isomorphic to the group algebra
$\ka[\Lambda]$. \esit

\subsection{Translation invariant function spaces}
\bsit\label{ftra} {\em Some characteristic-free function spaces.}
Consider the space $\mathcal F=\mathcal F(\Lambda,\ka)$ of all
$\ka$-valued functions on $\Lambda$. This is a
$\ka[\Lambda]$-module with respect to the natural action
$\ka[\Lambda]\times\mathcal F\to\mathcal F$, $(a,f)\longmapsto
a*f=\Delta_a(f)$.  A subspace $E\subseteq\mathcal F$ is a
$\ka[\Lambda]$-submodule if and only if it is translation
invariant i.e., $\tau_v (E)=E$ for all $v\in\Lambda$. A
one-dimensional subspace of $\cF$ is a $\ka[\Lambda]$-submodule if
and only if it is spanned by a character $\chi\in\Lambda^\vee$.
The following  function spaces are translation invariant. They
arise naturally in the framework of the Floquet-Bloch Theory, cf.\
\cite{Ku1, Ku2}.\esit

\bexas\label{exas}
\begin{enumerate}\item[(i)] The subspace $\ka[\Lambda]\subseteq\mathcal F$
of all functions with finite support.

\item[(ii)] Given a sublattice $\Lambda'\subseteq\Lambda$, we let
$\mathcal F_{\Lambda'}$ denote
the subspace of all $\Lambda'$-periodic functions i.e.,
$$\mathcal F_{\Lambda'}
=\{f\in\mathcal F\,
|\,\Lambda(f)\supseteq \Lambda'\}\,,$$ where
$$\Lambda(f)=\{v\in\Lambda\, |\,\tau_v(f)=f\}\,$$ stands for the
period lattice of $f$.

\item[(iii)] Given an isomorphism $\Lambda\cong \Z^s$ we define the subspace
$E_1\subseteq \mathcal F$ of all traces of polynomial functions
i.e.,
$$E_1=\{q|\Z^s\in\mathcal F\, |\, q\in\ka[x_1,\ldots,x_s]\}\,.$$
\item[(iv)] The subspace $E_2\subseteq \mathcal F$ of all pluri-periodic
functions,
$$E_2=\bigcup_{\Lambda'\subseteq
\Lambda,\,\ind_\Lambda\Lambda'<\infty} \mathcal F_{\Lambda'}\,,$$
is clearly translation invariant.

\item[(v)] The subspace $E_3\subseteq \mathcal F$ of all polynomially periodic
functions. Given an isomorphism $\Lambda\cong \Z^s$, a function
$g\in\mathcal F$ is called {\em polynomially periodic} if there
exists a sublattice $\Lambda'\subseteq \Lambda$ of finite index
and a set of polynomials $q_{\bar v}\in\ka[x_1,\ldots,x_s]$, where
$\bar v\in \Lambda/\Lambda'$, such that $$g(v)=q_{\bar v}(v)\quad
\mbox{if}\quad v\in \bar v\,.$$

\item[(vi)] The subspace $E_4\subseteq \mathcal F$ of all {\em Floquet functions}:
$$E_4=\{f=\sum_{i=1}^n\chi_i g_i\in\mathcal
F\, |\,\chi_i\in \Lambda^\vee=\Char (\Lambda,\ka^\times),
\,\,g_i\in E_3\}\,.$$

\item[(vii)] The subspace $E_5\subseteq \mathcal F$ spanned
by all $\Lambda'$-automorphic functions
on $\Lambda$, where $\Lambda'$ runs over the set of all finite
index sublattices of $\Lambda$. Given a sublattice
$\Lambda'\subseteq\Lambda$, a function $f\in\cF$ is called {\it
$\Lambda'$-automorphic} if,  for some character $\chi$ of
$\Lambda'$, $$(\tau_u f)(v)=\chi(u)f(v)\quad \forall u\in
\Lambda', \,\,\,\forall v\in \Lambda\,.$$

\item[(viii)] The subspace $E_6\subseteq \mathcal F$
of all finite-valued functions on $\Lambda$, etc.
\end{enumerate}\eexas

Due to the following proposition,  there are just three different
spaces among the function spaces $E_1,\ldots,E_6$. In particular,
in positive characteristic the Floquet functions are the same as
the pluri-periodic functions.

\bprop\label{pre} We have $$E_1\subsetneqq
E_2=E_3=E_4=E_5\subsetneqq E_6\,.$$\eprop

\bproof Letting $p=\Char (\ka)$ we observe that
$(x+p\lambda)^n=x^n$ for every monomial $x^n=\prod_{i=1}^s
x_i^{n_i}$ (regarded as a function in $x\in\Z^s=\Lambda$) and for
every $\lambda\in\Z^s$. Hence every polynomial function on
$\Lambda$ is $\Lambda'$-periodic, where $\Lambda'=p\Lambda$.
Therefore
$$E_1 \subseteq \mathcal F_{\Lambda'}
\subseteq E_2\,.$$ By a similar reason, $E_3,E_4 \subseteq E_2$.
For instance, if $g\in E_3$ is as in (v) then $p\Lambda'\subseteq
\Lambda(g)$. Indeed, for every $v\in\Lambda$ and $v'\in \Lambda'$
we have $v+pv'+\Lambda'=v+\Lambda'$. Hence
$$g(v+pv')=q_{\overline {v+pv'}}(v+pv')=q_{\overline {v+pv'}}(v)
=q_{\bar {v}}(v)=g(v)\,.$$ Thus $E_3\subseteq E_2$.

It is easily seen that the sum and the product of pluri-periodic
functions are again pluri-periodic. Since every character
$\chi\in\Lambda^\vee$ is $\Lambda_0$-periodic, where the
sublattice $\Lambda_0=\ker(\chi)$ is of finite index in $\Lambda$,
the inclusion $E_4\subseteq E_2$ holds. The inclusions
$E_2\subseteq E_3\subseteq E_4$ are evident. It follows that
$E_1\subseteq E_2=E_3=E_4\subseteq E_6$.

Every pluri-periodic function $f\in E_2$ (in particular, every
character of $\Lambda$) is $\Lambda(f)$-automorphic. Hence
$E_2\subseteq E_5$. Let furthermore $f\in E_5$ be a
$\Lambda'$-automorphic function as in (viii), where
$\Lambda'\subseteq\Lambda$ is a finite index sublattice. Since the
period lattice $\Lambda(f)$ contains the finite index sublattice
$\Lambda_0=\ker (\chi)\subseteq\Lambda'\subseteq\Lambda$, $f$ is
pluri-periodic. Thus $E_2=E_5$ and so $E_1\subseteq
E_2=E_3=E_4=E_5\subseteq E_6$.

It remains to show that $E_1\neq E_2\neq E_6$. The latter
assertion is easy and we leave it to the reader. As for the
former, it is enough to consider the case of a lattice
$\Lambda\cong\Z$ of rank $1$. In this case,  whenever $p\ge 3$, a
non-constant $2$-periodic function $f\in E_1$ cannot be a trace of
a polynomial. In other words, a 2-periodic trace of a polynomial
$q\in\ka [x]$ on $\Z$ must be constant. Indeed, since
$q(1)=q(1-p)=q(0)$ we have $q|2\Z=q|(2\Z+1)$. In the remaining
case $p=2$, a $3$-periodic function $f$ on $\Z$ with $f(0)=f(1)=0$
and $f(2)=1$ cannot be a trace of a polynomial over $\ka$, since
$f(0)\neq f(2)$. Thus indeed $E_1\neq E_2$.\eproof

\bsit\label{fftt} {\em The function space $\cF_p(\Lambda,\ka)$.}
In positive characteristic, the function space $E_2=E_3=E_4=E_5$
of all pluri-periodic functions on $\Lambda$ is not quite
appropriate for defining the Fourier transform. More suitable is
the translation invariant space
$$\cF_p(\Lambda,\ka)=\Span (\chi\vert \chi\in
\Lambda^\vee)\subseteq\cF$$ spanned by all characters of $\Lambda$
with values in $\ka^\times$. Of course, it depends  on the
characteristic $p$. In Proposition \ref{prop1} below we
characterize this subspace in terms of period lattices. To this
point we need the following definitions. \esit

\bdefi\label{satu} Let $\Lambda'\subseteq\Lambda$ be a sublattice
of finite index. We say that $\Lambda'$ is {\em $p$-saturated} if
$pv\in\Lambda'$ implies that $ v\in\Lambda'$. This holds if and
only if the index $\ind_\Lambda(\Lambda')$ is coprime to $p$, or,
equivalently, the quotient group $G=\Lambda/\Lambda'$ of order
$\ind_\Lambda(\Lambda')$ has no $p$-torsion. Clearly, the
intersection of two $p$-saturated sublattices is again
$p$-saturated. If $\Lambda'$ is $p$-saturated then every bigger
sublattice $\Lambda''\supseteq\Lambda'$ is, since
$\ind_\Lambda\Lambda''|\ind_\Lambda\Lambda'$.

If $\Lambda'$ is $p$-saturated then the dual group $G^\vee$ of $G$
is isomorphic to the group of characters $\Char (G,\ka^\times)$.
Otherwise $\Char (G,\ka^\times)\cong G^\vee/G^\vee(p)$, where
$G^\vee(p)$ is the Sylow $p$-component of $G^\vee$ (cf.\ \cite[\S
1]{Za2}).\edefi

\bsit\label{prep} Every character $\chi\in\Lambda^\vee$ regarded
as a function in $\cF_p (\Lambda,\ka)$ is pluri-periodic, with the
period lattice $\Lambda(\chi)=\ker (\chi)$. The quotient group
$$G(\chi)=\Lambda/\Lambda(\chi)\cong\chi(\Lambda)\subseteq
\ka^\times$$ is a finite cyclic group of order $\ord (\chi)$
coprime to $p$. Hence the period lattice $\Lambda(\chi)$ is
$p$-saturated (see \ref{satu}).

According to \ref{fftt} any function $f\in \cF_p (\Lambda,\ka)$ is
a linear combination of characters: $f=\sum_{i=1}^m
\alpha_i\chi_i$, where $\alpha_i\in\ka,\,\,\chi_i\in
\Lambda^\vee$. Therefore the lattice of periods of $f$,
$$\Lambda(f)\supseteq\bigcap_{i=1}^m \Lambda(\chi_i) $$
is also $p$-saturated.\esit

\bdefi\label{pp} Let $f\in\cF$ be a pluri-periodic function, that
is the period lattice $\Lambda(f)\subseteq\Lambda$ has finite
index (see Example \ref{exas}(iv)). We say that $f$ is {\em
$p$-pure} if $\Lambda(f)$ is $p$-saturated. We let $\cF_{pp}
(\Lambda,\ka)$ denote the space of all  $p$-pure pluri-periodic
functions on $\Lambda$ with values in $\ka$. This is indeed a
linear space because $\Lambda(f+g)\supseteq\Lambda(f)\cap
\Lambda(g)$ is $p$-saturated if $\Lambda(f)$ and $ \Lambda(g)$
are, see \ref{satu}. Furthermore, the space $\cF_{pp}
(\Lambda,\ka)$ is translation invariant, because
$\Lambda(\tau_v(f))=\Lambda(f)$ $\forall v\in\Lambda$. So this is
a convolution module over $\ka[\Lambda]$.\edefi

\brem\label{non-pure} We have shown in Example \ref{exas}  that
the trace  on $\Lambda$ of a polynomial $q\in\ka[x_1,\ldots,x_s]$
is $\Lambda'$-periodic, where $\Lambda'=p\Lambda=\langle
pv_1,\ldots,pv_s\rangle$. However, the sublattice
$\Lambda'\subseteq\Lambda$ of finite index is not $p$-saturated.
The function $q|\Lambda$ is $p$-pure if and only if this is a
constant function.\erem

\bprop\label{prop1} $\cF_{pp} (\Lambda,\ka)=\cF_p (\Lambda,\ka)$.
\eprop

\bproof The inclusion $\cF_p (\Lambda,\ka)\subseteq\cF_{pp}
(\Lambda,\ka)$ is immediate from \ref{prep} and \ref{pp}. It
remains to show the converse inclusion. Let $f$ be a
pluri-periodic function with period lattice
$\Lambda(f)\subseteq\Lambda$ of index coprime to $p$. Then $f=\bar
f\circ\pi$ is the pullback of a function $\bar f: G(f)\to\ka$,
where $G(f)=\Lambda/\Lambda(f)$ is a finite abelian group of order
coprime with $p$ and $\pi:\Lambda\to G(f)$ is the canonical
surjection. The function space $\mathcal F (G(f),\ka)$ being
spanned by characters of $G(f)$ with values in $\ka^\times$ (see
\cite[Lemma 2.1(a)]{Za2}) there is a presentation
\be\label{prst}\bar f=\sum_{i=1}^n
\alpha_i\overline\chi_i\quad\mbox{with}\quad
\overline\chi_i\in\Char (G(f),\ka^\times)\quad\mbox{and}\quad
\alpha_i\in\ka\,.\ee Thus $f=\sum_{i=1}^n \alpha_i\chi_i$, where
the characters $\chi_i=\overline\chi_i\circ\pi\in\Lambda^\vee$ are
the pullbacks of $\overline\chi_i$, $i=1,\ldots,n$. This yields
the desired conclusion. \eproof

\subsection{Fourier transform}
\bsit\label{futr} {\em The Fourier transform on $\ka[\Lambda]$}
sends the group ring $\ka[\Lambda]$ isomorphically onto the
coordinate ring of the torus $\mathcal O(\Lambda^\vee)$. Given a
base $\mathcal V$ of $\Lambda$ and the identifications
$\Lambda\cong\Z^s$ and $\Lambda^\vee\cong\T^s$ as above, $F$ is
given by
$$F:f\longmapsto \widehat f (z)=
\sum_{\lambda\in\Z^s} f(\lambda)z^{\lambda}\,.$$ It sends
$\ka[\Lambda]\cong\ka[\Z^s]$ isomorphically onto the algebra of
Laurent polynomials
$$\mathcal O(\Lambda^\vee)\cong\mathcal O(\T^s)=\ka [z,z^{-1}]=\ka
[z_1,\ldots,z_n,z_1^{-1},\ldots,z_n^{-1}]\,.$$ In particular, the
delta-function $\delta_v\in \ka[\Lambda]$ is sent to the Laurent
monomial $\widehat{\delta_v}=z^{\lambda(v)}\in \mathcal O(\T^s)$,
which is a character of $\Lambda$. Or, in coordinate-free form,
$\widehat{\delta_v}=\eval_v$ regarded as a function on the
characters of $\Lambda$. The equality
$\delta_v*\delta_w=\delta_{v+w}$ in $\ka[\Lambda]$ amounts to
$z^{\lambda(v)} z^{\lambda(w)}=z^{\lambda(v)+\lambda(w)}$ in
$\mathcal O(\T^s)$. The inverse Fourier transform $F^{-1}:\ka
[z,z^{-1}]\to \ka[\Z^s]$ sends a Laurent polynomial
$\sum_{\lambda\in\Z^s} f(\lambda)z^\lambda$ into its coefficient
function $f\in\ka[\Z^s]$. \esit

\bsit\label{futr1} {\em Fourier transform on $\cF_p
(\Lambda,\ka)$.} We let $\cF_{\rm fs} (\Lambda^\vee,\ka)$ denote
 the space of all functions on $\Lambda^\vee$ with finite
support. It comes equipped with the structure of $\mathcal
O(\Lambda^\vee)$-module.

Our next goal is to define a Fourier transform on the space $\cF_p
(\Lambda,\ka)\subseteq\cF$ so to obtain an isomorphism
$$F: \ka[\Lambda]\oplus\cF_p
(\Lambda,\ka)\stackrel{\cong}{\longrightarrow}\cO(\Lambda^\vee)\oplus\cF_{\rm
fs} (\Lambda^\vee,\ka)\,$$  of $\ka[\Lambda]$- and
$\cO(\Lambda^\vee)$-modules, respectively. The following
observations will be useful in the sequel.\esit

\bsit\label{dul} Suppose that $\Lambda'$ is  $p$-saturated, and
let $\pi:\Lambda\to G=\Lambda/\Lambda'$ be the canonical
surjection. The dual group $G^\vee\cong\Char (G,\ka^\times)$ of
the finite group $G=\Lambda/\Lambda'$ can be realized as the
finite subgroup
$${\Lambda'}^\bot=\{\chi\in\Lambda^\vee\,|\,\Lambda' \subseteq
\ker \chi\} \subseteq \Lambda^\vee\,$$ via the pullback
$$\pi^\vee:G^\vee\stackrel{\cong}\longrightarrow {\Lambda'}^\bot,\qquad
G^\vee\ni\bar\chi\longmapsto\chi=\pi^\vee(\bar\chi)=\bar\chi\circ\pi\in
{\Lambda'}^\bot\,.$$ Letting $f=\bar f\circ\pi\in
\cF_{\Lambda'}(\Lambda,\ka)$ and $\chi=\bar\chi\circ\pi\in
\Lambda'^\bot$ from (\ref{prst}) we obtain a presentation
\be\label{ftr1} f=\sum_{\chi\in
\Lambda'^\bot}\alpha(\chi)\cdot\chi \,.\ee Clearly, every function
$f\in \cF_{\Lambda'}(\Lambda,\ka)$ admits such a presentation
i.e., $\cF_{\Lambda'}(\Lambda,\ka)=\Span
(\chi\vert\chi\in\Lambda'^\bot)$. \esit

\bdefi\label{ftpp}  For a function $f\in\cF_p (\Lambda,\ka)$ we
define its Fourier transform as \be\label{ftr2} F: f=\sum_{i=1}^m
\alpha_i\chi_i\longmapsto \widehat f=\sum_{i=1}^m
\alpha_i\delta_{\chi_i^{-1}}\,.\ee Since $\supp (\widehat f)
\subseteq \Lambda(f)^\bot$ (see \ref{dul}) we have $\widehat
f\in\cF_{\rm fs} (\Lambda^\vee,\ka)$. Thus the Fourier transform
yields an isomorphism
$$F:\cF_p (\Lambda,\ka)\stackrel{\cong}{\longrightarrow}\cF_{\rm fs}
(\Lambda^\vee,\ka),\qquad \chi\longmapsto \delta_{\chi^{-1}}\,,$$
with inverse
$$\cF_{\rm fs} (\Lambda^\vee,\ka)\ni\varphi=\sum_{\chi\in\supp
(\varphi)} \varphi(\chi)\cdot\delta_\chi\longmapsto
F^{-1}(\varphi)=\sum_{\chi\in\supp
(\varphi)}\varphi(\chi)\cdot\chi^{-1}\in \cF_p
(\Lambda,\ka)\,.$$\edefi

\bprop\label{isomm} For every $a\in\ka[\Lambda]$ and for every
$f\in\cF_p(\Lambda,\ka)$ we have \be\label{b}
\widehat{a*f}=\widehat{a}\cdot \widehat{f}\,.\ee \eprop

\bproof The space $\ka[\Lambda]$ being spanned by delta functions
and $\cF_p(\Lambda,\ka)$ by characters of $\Lambda$, and the
Fourier transform being linear, it is enough to verify (\ref{b})
in the particular case where $a=\delta_v$ ($v\in\Lambda$) and
$f=\chi \in\Lambda^\vee$. We have
$$\widehat{\delta_v*\chi}=\widehat{\tau_{-v}(\chi)}
=\chi^{-1}(v)\cdot\widehat{\chi}=\widehat{\delta_v}(\chi^{-1})\cdot\delta_{\chi^{-1}}=
\widehat{\delta_v}\cdot\delta_{\chi^{-1}}=\widehat{\delta_v}\cdot\widehat{\chi}\,.$$
Indeed, $\widehat{\delta_v}=\eval_v\in\cF(\Lambda^\vee,\ka)\,$. In
particular
$$\chi^{-1}(v)=\eval_v(\chi^{-1})=\widehat{\delta_v}(\chi^{-1})\,.$$
Now (\ref{b}) follows.  \eproof

\section{Convolution operators: the scalar case}
\subsection{$\Delta_a$-harmonic functions}
\bsit\label{aharm} Let $a\in\ka[\Lambda]$. A function $f\in\cF$
satisfying the equation
$$\Delta_a(f)=a*f=0$$ will be called {\it $\Delta_a$-harmonic}.
Choosing a base of $\Lambda$ we identify $\Lambda$ with $\Z^s$ and
$\Lambda^\vee$ with the torus $\T^s$. With every convolution
operator $\Delta_a$ on $\Lambda=\Z^s$ we associate its {\em
symbolic variety} $\Sigma_a=\Sigma_{\Delta_a}$. This is an affine
hypersurface in the dual torus $\Lambda^\vee\cong\T^s$ defined as
$$\Sigma_a=\{z\in\T^s\,|\,\widehat a(z^{-1})=0\}\,,$$ where  $\widehat a=F(a)$
is a Laurent polynomial in $z=(z_1,\ldots,z_s)$, see
\cite[2.7]{Za2}. Every character $\chi\in\Sigma_a$ viewed as a
function on $\Lambda=\Z^s$ is $\Delta_a$-harmonic. Indeed, for
$z\in\Sigma_a$ we have
$$\Delta_a(z^\lambda)=\sum_{v\in\Lambda} a(v)\tau_{-v} (z^\lambda)=
\sum_{v\in\Lambda}a(v)z^{\lambda-v}=z^\lambda\widehat a
(z^{-1})=0\,.$$  \esit

The following proposition is an analog of Theorems 1-3 in
\cite{Ku2} in the case of scalar  functions.

\bprop\label{scca}\begin{enumerate}\item[(a)] If $a\neq 0$ then a
nonzero $\Delta_a$-harmonic function on $\Lambda$ cannot have
finite support.
 \item[(b)] Every $\Delta_a$-harmonic function
$f\in\cF_p(\Lambda,\ka)$ is a linear combination of
$\Delta_a$-harmonic characters $\chi\in\Sigma_a$. In other words,
$$\ker (\Delta_a|\cF_p(\Lambda,\ka))=\Span
(\chi\,|\,\chi\in\Sigma_a)\,.$$ \end{enumerate} \eprop

\bproof To show (a) suppose that $f\in\ka[\Lambda]$ is
$\Delta_a$-harmonic. Then by Proposition \ref{isomm} ${\widehat
a}\cdot {\widehat f}=0$ and so ${\widehat a}, \,{\widehat
f}\in\ka[z,z^{-1}]$ are zero divisors. However, the algebra
$\ka[z,z^{-1}]$ of Laurent polynomials being an integral domain
this implies that $f=0$. Indeed, ${\widehat a}\neq 0$ by our
assumption.

(b) If $a=c\delta_v$ i.e., $\Delta_a$ is proportional to a shift
then the only $\Delta_a$-harmonic function in $\cF$ is the zero
function. Assume further that $a\neq c\delta_v$ i.e., $\widehat a$
is not a Laurent monomial and so $\Sigma_a\neq\emptyset$. Let
$f\in\cF_{\Lambda'}$ be $\Delta_a$-harmonic and
$\Lambda'$-periodic, where $\Lambda'\subseteq\Lambda$ is a
sublattice of a finite index coprime to $p$, see Proposition
\ref{prop1}. Then $f$ is the pullback of a function $\bar f:
G=\Lambda/\Lambda'\to\ka$. We let $a_*:G\to\ka$ denote the
pushforward  function defined by
$$ a_*(v+\Lambda')=\sum_{v'\in \Lambda'} a(v+v')\,.$$ This is a unique
function in $\cF(G,K)$ satisfying the condition $a_*\circ
\pi=a*\delta_{\Lambda'}$, where $\pi: \Lambda\twoheadrightarrow G$
is the canonical surjection and $\delta_{\Lambda'}$ stands for the
characteristic function of $\Lambda'\subseteq\Lambda$ (see
\cite[(7)]{Za2}). Since $f$ is $\Delta_a$-harmonic, by the
Pushforward Lemma 1.1 in \cite{Za2}, $\Delta_{a_*} (\bar f)=0$. It
follows that $\bar f$ is a linear combination of $a_*$-harmonic
characters of $G$, see \cite[Corollary 2.2]{Za2}. Their pullbacks
are $\Delta_a$-harmonic characters of $\Lambda$, and so $f$ is a
linear combination of the latter characters, as stated in (b).
\eproof

There are the easy equivalencies:
$$a=0\quad\Longleftrightarrow\quad\Sigma_a=\T^s\quad\Longleftrightarrow\quad
\ker (\Delta_a|\cF_p(\Lambda,\ka))=\cF_p(\Lambda,\ka)\,.$$ In the
opposite direction, we have the following result (cf.\
\cite[Corollary 2]{Ku1}).

\bcor\label{co1} There is no nonzero $\Delta_a$-harmonic function
on $\Lambda$ if and only if $\Sigma_a=\emptyset$, if and only if
$\Delta_a$ is proportional to a shift, if and only if $\Delta_a$
is invertible.\ecor

\subsection{The Floquet-Fermi
hypersurfaces} This notion is borrowed in \cite[Definition
2]{KuVa}.

\bdefi\label{fs} For a convolution operator $\Delta_a$, where
$a\in\ka[\Lambda]$, its {\em Floquet-Fermi hypersurface} of level
$\mu\in\ka$ is
$$\Sigma_{a,\mu}=\{z\in\T^s\,|\,\widehat a(z^{-1})=\mu\}\,.$$
Actually $\Sigma_{a,\mu}=\Sigma_{a-\mu\delta_0}$ is the symbolic
hypersurface of the operator $\Delta_a-\mu\Id$. We let
$$E_{a,\mu}=\Span (\chi\,|\,\chi\in\Sigma_{a,\mu})\,$$
denote the eigenspace of $\Delta_a|\cF_p(\Lambda,\ka)$ with
eigenvalue $\mu$.
 \edefi

The following proposition is straightforward. Nevertheless, we
provide a short argument.

\bprop\label{sp1} We have
$$\cF_p(\Lambda,\ka)=\bigoplus_{\mu\in\ka} E_{a,\mu}\,.$$
In particular, the convolution operator $\Delta_a$ on the space
$\cF_p(\Lambda,\ka)$ is diagonalizable in the basis of characters.
\eprop \bproof Choosing a base $(v_1,\ldots,v_s)$ of the lattice
$\Lambda$, we identify $\Lambda$ with $\Z^s$ and $\Lambda^\vee$
with the torus $\T^s$. Letting $\tau_i$ denote the shift by
$-v_i$, the convolution operator $\Delta_a$ can be written as a
Laurent polynomial in the basic shifts $\tau_1,\ldots,\tau_s$.
Namely, $$\Delta_a=\sum_{\lambda\in\Z^s}
a(\lambda)\tau_{-\lambda}=\sum_{\lambda=(\lambda_1,\ldots,\lambda_s)\in\Z^s}
a(\lambda)\prod_{i=1}^s\tau_i^{\lambda_i}=\widehat
{a}(\tau_1,\ldots,\tau_s)\,.$$ Using the equalities
$\tau_i(z^\lambda)=z_i^{-1}\cdot z^\lambda$, for every
$z\in\ka^\times$ we obtain
$$\Delta_a (z^\lambda)=\widehat
{a}(\tau_1,\ldots,\tau_s)(z^\lambda)=\widehat {a}(z^{-1})\cdot
z^\lambda\,.$$ Now both assertions follow. \eproof

\brems\label{muop} 1. Here is an alternative approach. According
to Proposition \ref{isomm} the Fourier transform $F$ sends the
convolution operator $\Delta_a|\cF_p(\Lambda,\ka)$ to the operator
of multiplication by the function $\widehat a$ on the space
$\cF_{\rm fs}(\Lambda^\vee,\ka)=F(\cF_p(\Lambda,\ka))$. The latter
operator being semi-simple, $\Delta_a$ is semi-simple too.

2. Let us  identify $\Lambda$ with $\Z^s$ and $\Lambda^\vee$ with
the torus $\T^s$.  Notice that the spectrum of
$\Delta_a|\cF_p(\Lambda,\ka)$ coincides with the range of the
Laurent polynomial $$\widehat a(z^{-1})=\sum_{\lambda\in\Z^n}
a(\lambda)z^{-\lambda}\,.$$ Hence \bnum\item[$\bullet$] $\spec
(\Delta_a) = \{\mu\}$ if and only if
$a=\mu\delta_0$;\item[$\bullet$] $\spec (\Delta_a) = \ka^\times$
if and only if $a=\alpha \delta_v$ for some $\alpha\in\ka^\times$
and some $v\in\Lambda\setminus\{0\}$; \item[$\bullet$] $\spec
(\Delta_a) =\ka$ otherwise.\enum Furthermore, if $a\neq \alpha
\delta_v$ then \be\ba{ll} E_{a,\mu}=\Span
(z^{-\lambda}\,|\,\widehat
a(z)=\mu)=\Big\{f\in\cF_p(\Lambda,\ka)\,|\,
f(\lambda)=\sum_{z\in\T^s} \alpha(z)z^{-\lambda}\,:
\\
\hfill \alpha\in\cF_{\rm fs}(\T^s,\ka)\,\,\mbox{
s.t.}\,\,\supp(\alpha)\subseteq\widehat a^{-1}(\mu)\Big\}\,.
\\
\ea\ee If $s=\rank(\Lambda)>1$ then every level set $\widehat
a^{-1}(\mu)$, where $\mu\in\spec (\Delta_a)$, is an infinite
countable set (indeed, the algebraically closed field
$\ka=\overline{\F_p}$ is infinite countable). Hence for every
$\mu\in\ka$ the eigenspace $E_{a,\mu}$ has a countable basis of
characters; in particular $\dim(E_{a,\mu})=+\infty$. While for
$s=1$ all these eigenspaces are finite dimensional. \erems

\section{Convolution operators: the vector case}

Letting $\A^n=\A^n_\ka$ be the affine $n$-space over $\ka$ and
 $(e_1,\ldots,e_n)$ be the  canonical basis of $\A^n$,
we consider the space $\cF^n=\cF(\Lambda,\A^n)$ of all
$\ka$-valued vector functions on $\Lambda$, and its subspaces
$\ka[\Lambda]^n$ of all vector functions on $\Lambda$ with finite
support and $\cF_p(\Lambda,\A^n)$ of all  $p$-pure pluri-periodic
vector functions on $\Lambda$. The vector functions on
$\Lambda=\Z^s$ of the form $f(\lambda)=\chi\cdot u=z^\lambda\cdot
u$, where $z\in\T^s=\Lambda^\vee$ and $u\in\A^n$, will be called
{\em elementary}. These are  analogs  in positive characteristic
of the classical Bloch functions. Thus any $f\in
\cF_p(\Lambda,\A^n)$ is a linear combination of elementary
functions.

\bsit\label{newl0} Identifying as before $\Lambda$ with $\Z^s$ and
$\Lambda^\vee$ with $\T^s$, the Fourier transform sends a vector
function $f\in \ka[\Lambda]^n$ into the vector Laurent polynomial
$$\widehat f(z)=\sum_{\lambda\in\Z^n}
f(\lambda)z^{\lambda},\qquad z\in\T^s\,,$$ where
$f(\lambda)\in\A^n$ $\forall \lambda\in\Z^s$. In particular, a
basic vector delta-function $\delta_{i,\lambda}=\delta_\lambda
\cdot e_i\in \ka[\Lambda]^n$, where $\lambda\in\Lambda=\Z^s$, is
sent to the vector Laurent monomial $\widehat
{\delta_{i,\lambda}}=z^{\lambda}\cdot e_i$ ($i=1,\ldots,n$). A
pluri-periodic vector function
$$f=\sum_{i=1}^m \alpha_i\chi_i\in \cF_p(\Lambda,\A^n)\,,$$
where $\alpha_i\in\A^n$ and $\chi_i\in\Lambda^\vee$, is sent to
the vector function on $\Lambda^\vee$ with finite support
$$\widehat f=\sum_{i=1}^m \alpha_i\delta_{\chi_i^{-1}}\in \cF_{\rm
fs}(\Lambda^\vee,\A^n)\,.$$
 \esit

\subsection{$\Delta$-harmonic vector functions}\label{newl1}
Consider an endomorphism $\Delta\in\End\,(\ka[\Lambda]^n)$
commuting with shifts on $\Lambda$. Letting
$$\Delta(\delta_{i,0})=
a_i=(a_{i1},\ldots,a_{in})\in\ka[\Lambda]^n,\quad
i=1,\ldots,n\,,$$ we get a square matrix $A=(a_{ij})$ of order $n$
with entries in $\ka[\Lambda]$. The images of the remaining basic
vector functions $\delta_{i,v}$ can be recovered as suitable
shifts of the vectors $a_i\,$:
$$\Delta (\delta_{i,v})=\Delta
(\tau_{-v}\delta_{i,0})=\tau_{-v} (\Delta(\delta_{i,0}))=\tau_{-v}
(a_i)\,.$$ If $f=(f_1,\ldots,f_n)\in \ka[\Lambda]^n$ and
$g=\Delta(f)=(g_1,\ldots,g_n)\in\ka[\Lambda]^n$ then $g=A*f$ i.e.,
\be\label{new} g_i=\sum_{i=1}^n a_{ij}*f_j,\qquad
i=1,\ldots,n\,.\ee In other words, $\Delta$ acts on
$\ka[\Lambda]^n$ as a square matrix of convolution operators
$\Delta_A=(\Delta_{a_{ij}})$. Such operator $\Delta$ can be
applied to an arbitrary vector function $f=(f_1,\ldots,f_n)\in
\cF(\Lambda,\A^n)=\cF^n$. This yields an extension of $\Delta$ to
$\cF^n$, which we denote by the same letter. The extended
endomorphism $\Delta\in\End(\cF^n)$ leaves invariant any
translation invariant subspace. In particular, it acts on the
subspace $\cF_p(\Lambda,\A^n)\subseteq\cF^n$.

The Fourier transform sends $A$ into a matrix of Laurent
polynomials $\widehat A=(\widehat a_{ij})$. By Proposition
\ref{isomm} the relation $g=A*f$, where $f,g\in
\cF_{p}(\Lambda,\A^n)$, is transformed into
$$\widehat g=\widehat A\cdot \widehat f\,,$$ where
$\widehat f,\widehat g\in \cF_{\rm fs}(\Lambda^\vee,\A^n)$.

\medskip

The following is an analog of Theorems 1-3 in \cite{Ku1} for
vector functions, with a similar proof  (cf.\ also Theorem 8 in
\cite{Ku2}).

\bthm\label{kuch} Let $\Delta\in\End\,(\ka[\Lambda]^n)$ be an
endomorphism commuting with shifts.  Consider the equation
\be\label{eq0} \Delta(f)=0\,,\ee where $\Delta$ acts on $\cF^n$ as
in (\ref{new}). Then the following hold.
\begin{enumerate}\item[(a)] There exists a non-zero solution of (\ref{eq0})
with finite support $f\in \ka[\Lambda]^n$  if and only if
$\det\widehat A=0$.
\item[(b)] Let $\Lambda'\subseteq\Lambda$ be a sublattice of finite
index coprime to $p$. There exists a $\Lambda'$-periodic solution
$f\in\cF_{\Lambda'}(\Lambda,\A^n)$ of (\ref{eq0}) if and only if
there exists a character $\chi\in \Lambda'^\bot$ with $\det
\widehat A (\chi^{-1})=0$. Every such solution is a linear
combination of elementary solutions \be\label{vf} f=\sum_{i=1}^m
\chi_i \cdot u_i\,,\ee where $\chi_i\in\Lambda'^\bot$ with $\det
\widehat A (\chi_i^{-1})=0$ and $u_i\in \ker (\widehat A
(\chi_i^{-1}))$  $\forall i=1,\ldots,m$. Vice versa, any vector
function as in (\ref{vf}) is a $\Lambda'$-periodic solution of
(\ref{eq0}).
\end{enumerate}\ethm

\bproof (a) Choosing a base of $\Lambda$ we identify $\Lambda$
with $\Z^s$ and $\Lambda^\vee$ with $\T^s$. Applying the Fourier
transform to (\ref{eq0}) yields the equation \be\label{eq}\widehat
A\cdot \widehat f=0\,.\ee We can treat (\ref{eq}) as a linear
system over the rational function field $\ka(z)$. By the usual
linear algebra argument, this system admits a rational solution
$\widehat f\in [\ka(z)]^n=[\ka(z_1,\ldots,z_s)]^n$ if and only if
$\det\widehat A=0$. Multiplying such a solution $\widehat f$ by a
suitable polynomial $q\in\ka[z]$ we obtain a polynomial solution
of (\ref{eq}), hence a solution in vector Laurent polynomials. In
turn, (\ref{eq0}) admits a solution with finite support.

(b) If ${f}=({f_1},\ldots, {f_n})$ is a nonzero
$\Lambda'$-periodic solution of (\ref{eq0}) then $\widehat
{f}=(\widehat {f_1},\ldots, \widehat {f_n})$ is a nonzero solution
of (\ref{eq}) supported on ${\Lambda'}^\bot$. Hence $\widehat f$
does not vanish at some point $\chi\in {\Lambda'}^\bot$ and so
$(\det \widehat A)(\chi)=0$. However, any vector $u\in \ker
(\widehat A (\chi^{-1}))$ yields an elementary solution
$f=\chi\cdot u$ of (\ref{eq0}). Now the first assertion in (b)
follows. The last assertion of (b) is evident.

To show the second assertion in (b), let us consider the finite
abelian group $G=\Lambda/\Lambda'$ of order coprime to $p$. Every
$\Lambda'$-periodic vector function
$f\in\cF_{\Lambda'}(\Lambda,\A^n)$ is pullback of a vector
function $\bar f\in \cF(G,\A^n)$. The image $g=A*f$ being also
$\Lambda'$-periodic, it is pullback of a function $\bar g\in
\cF(G,\A^n)$. The endomorphism $\Delta_*:\bar f\longmapsto \bar g$
is given by the pushforward matrix $A_*=(a_{ij*})$, where
$a_{ij*}\in \cF(G,\ka)$ $\forall (i,j)$. The latter acts via
$$\bar g_i=\sum_{j=1}^n a_{ij*}*\bar f_j,\qquad i=1,\ldots,n\,$$
(cf.\ the Pushforward Lemma 1.1 in \cite{Za2}). Since $p\nmid
\ord\,G$ the Fourier transform is well defined on the space
$\cF(G,\ka)$ (cf.\ \cite{Za2}). Applying the Fourier transform on
$G$ yields a matrix function $\widehat {A_*}=(\widehat {a_{ij*}})$
on the dual group $G^\vee$ which acts on vector functions on
$G^\vee$ via the usual matrix multiplication. Thus (\ref{eq0})
leads to equation \be\label{eq1} \widehat {A_*}\cdot \widehat
{\bar f}=0\,.\ee Any solution $f$ of (\ref{eq0})  gives a solution
$\widehat {\bar f}$ of (\ref{eq1}). Writing this solution as a
linear combination of delta-functions yields a presentation of
$\bar f$ as a linear combination of elementary solutions on $G$.
The pullback of this presentation gives (\ref{vf}). \eproof

\bdefi\label{mult} Consider  as before an endomorphism
$\Delta\in\End(\ka[\Lambda]^n)$ commuting with shifts, where
$n\in\N$. Similarly to the scalar case (where $n=1$) we call the
hypersurface
$$\Sigma_\Delta=\{z\in\Lambda^\vee\,|\,\det \widehat A(z^{-1})=0\}$$
the {\em symbolic variety} of $\Delta$ (cf.\ \ref{aharm}). A point
$z\in\Sigma_\Delta$ corresponding to a character $\chi=z^\lambda$
of $\Lambda\cong\Z^s$ is called a {\it multiplier} of (\ref{eq0}).
In case (b) of Theorem \ref{kuch} the multipliers run over the
intersection $\Sigma_\Delta\cap {\Lambda'}^\bot$. Thus every
$\Lambda'$-periodic solution $f$ of (\ref{eq0}) as in (\ref{vf})
can be written as
$$f(\lambda)=\sum_{z\in \Sigma_\Delta\cap {\Lambda'}^\bot}
z^\lambda\cdot u(z)\,.$$ \edefi

\bcor\label{3co} Under the assumptions of Theorem \ref{kuch} the
equation (\ref{eq0}) admits a nonzero solution
$f\in\cF_p(\Lambda,\ka)$ if and only if it admits such an
elementary solution, if and only if $\Sigma_\Delta\neq\emptyset$.
Moreover, $$\ker (\Delta)=\Span \left(z^\lambda\cdot
u\,|\,z\in\Sigma_\Delta, u\in\ker \widehat
A(z^{-1})\right)\,.$$\ecor

\brem\label{core} The finite subgroups ${\Lambda'}^\bot$ exhaust
the  dual torus $\Lambda^\vee\cong\T^s$ when $\Lambda'$ runs over
all sublattices of $\Lambda$ of finite indices coprime to $p$.
Given $\Delta\in\End \ka[\Lambda]^n$ commuting with shifts, it
would be interesting to find the counting function $\card
(\Sigma_\Delta\cap {\Lambda'}^\bot)$ as a suitable analog of the
Weil zeta-function. \erem

\subsection{The Jordan form of a matrix convolution operator}\label{rjf}
Over an algebraically closed field,
every square
matrix  admits a Jordan basis \cite[Ch. VII, \S 5.2, Prop.
5]{Bou}. We construct below such a basis for any matrix
convolution operator.

\bdefi\label{jd} Given an operator $\Delta\in\End(\ka[\Lambda])$
commuting with shifts,  for every $\mu\in\ka$ we consider its {\em
Floquet-Fermi hypersurface}
$$\Sigma_{\Delta,\mu}=\{z\in\T^s\,|\,\det(\widehat
A(z^{-1})-\mu\cdot I_n)=0\}\,,$$ the corresponding eigenspace
$$E_{\Delta,\mu}=\ker ((\Delta-\mu\cdot \Id)|\cF_p(\Lambda,\ka))\,,$$
and the generalized  eigenspace
$$E_{\Delta}^{(\mu)}=\ker ((\Delta-\mu\cdot\Id)^n|\cF_p(\Lambda,\ka))\,.$$
\edefi

The proof of the following theorem is straightforward.

\bthm\label{spde}\begin{enumerate} Let
$\Delta\in\End(\ka[\Lambda])$ be an endomorphism commuting with
shifts, extended  via (\ref{new}) to the space
$\cF_p(\Lambda,\ka)$ as a matrix convolution operator. Then the
following hold.
\item[(a)] There is a decomposition
$$\cF_{p}(\Lambda,\A^n)=\bigoplus_{\mu\in\ka} E_{\Delta}^{(\mu)}\,.$$
\item[(b)] Every generalized eigenspace
$E_{\Delta}^{(\mu)}\subseteq\cF_p(\Lambda,\ka)$
is generated by elementary vector-functions. Namely,
$$E_{\Delta}^{(\mu)}=\Span \left(z^\lambda\cdot u \,|\,z\in
\Sigma_{\Delta,\mu},\,\,u \in\ker(A(z^{-1})-\mu\cdot
I_n)^n\right)$$
$$=\Big\{ f\in \cF_{p}(\Lambda,\A^n)\,|\,f(\lambda)=
\sum_{z\in\Sigma_{a,\mu}}
\alpha(z)z^\lambda,\quad\mbox{where}\quad\alpha\in\cF_{\rm
fs}(\Lambda^\vee,\A^n)\quad\mbox{is s.t.}$$$$\alpha(z)\in\ker
(\widehat A(z^{-1})-\mu\cdot I_n)^n\,\,\,\forall z\in \supp
(\alpha)\Big\}\,.$$
\item[(c)] The induced endomorphism $\Delta\in \End(\cF_p(\Lambda,\ka))$
admits a Jordan basis consisting of elementary vector functions.
More precisely, fixing for every $z\in\T^n$ a Jordan basis
$(u_1(z),\ldots,u_n(z))$ for $\widehat A(z^{-1})$  in $\A^n$ we
obtain a Jordan basis $(z^\lambda\cdot u_i(z)\,|\,z\in\T^n)$ for
$\Delta$ in $\cF_p(\Lambda,\A^n)$.
\end{enumerate}\ethm

\brem\label{coco} \`A priori, the choice of the vector functions
$u_i(z)$ in (c) does not suppose any specific dependence on
$z=(z_1,\ldots,z_s)$. However, they can be chosen as
(non-rational, in general) algebraic vector functions in $z$.
\erem

\subsection{Convolution equations over finite fields}\label{gfi}
Taking suitable traces we can obtain solutions of the linear
equation (\ref{eq0}) with values in a finite field. We discuss
below this matter in more detail.

Given a square matrix function $A=(a_{ij})$ of order $n$ with
entries $a_{ij}\in\ka[\Lambda]$, we consider as before the matrix
convolution equation \be\label{mace} A*f=0,\quad\mbox{where}\quad
f=(f_1,\ldots,f_n)\in\cF_p(\Lambda,\A^n)\,.\ee The proof of the
following proposition follows the lines of the proof of Theorem
4.4 in \cite{Za2}.

\bprop\label{mceff} Consider a finite Galois field
$GF(q)=\F_q\subset\ka$, where $q=p^\sigma$. Suppose that
$a_{ij}\in GF(q)[\Lambda]$ $\forall i,j=1,\ldots,n$. If
(\ref{mace}) admits a nonzero solution $f\in\cF_p(\Lambda,\A^n)$,
then it admits also a nonzero solution $\tilde f\in
\cF_p(\Lambda,\A^n)$ with values in the field $GF(q)$. \eprop

\bproof The equation (\ref{mace}) is equivalent to the system
$$\sum_{j=1}^n a_{ij}*f_j=0,\qquad i=1,\ldots,n\,.$$ Since $GF(q)$
is the subfield of invariants of the Frobenius automorphism
$z\longmapsto z^q$ acting on $\ka$, applying this automorphism
yields the equalities $$\sum_{j=1}^n a_{ij}*f^q_j=0,\qquad
i=1,\ldots,n\,.$$ That is
$f^q=(f_1^q,\ldots,f_n^q)\in\cF_p(\Lambda,\A^n)$ is again a
solution of (\ref{mace}).

The function $f_j$ being pluri-periodic its image
$f_j(\Lambda)\subseteq\ka$ is finite. Let $q(f)=p^{r(f)}$ be such
that $\GF(q(f))$ is the smallest subfield of $\ka$ containing
$\GF(q)$ and all the images $f_j(\Lambda)$, $j=1,\ldots,n$. We let
$\Tr(f)=(\Tr(f_1),\ldots,\Tr(f_n))$, where
$$\Tr(f_j)=\Tr_{GF(q(f)):GF(q)}(f_j):
=f_j+f_j^q+\ldots+f_j^{q^{r(f)-1}}\in\cF_p(G,\ka)\,.$$ The trace
$\Tr(f)$ is fixed by the Frobenius automorphism, hence it
represents a solution of (\ref{mace}) with values in the finite
field $GF(q)$. It remains to show that such a solution can be
chosen nonzero.

Let $\Delta\in\End (\cF_p(\Lambda,\A^n))$ denote the matrix
convolution operator $\Delta (f)=A*f$. By Corollary \ref{3co}
$\ker (\Delta)$ is spanned by elementary solutions
$f=z^\lambda\cdot u$, where $z\in\Sigma_\Delta\subseteq\T^s$. By
our assumption $\ker (\Delta)\neq (0)$, hence $\Sigma_\Delta\neq
\emptyset$. Choosing $z\in \Sigma_\Delta$ we let
$$r(z)=\min\left\{r\,|\,z^\lambda,
(z^q)^\lambda,\ldots,(z^{q^{r-1}})^\lambda\quad\mbox{are all
distinct}\right\}\,.$$ The entries of the matrix $\widehat
A(z^{-1})$ being Laurent polynomials with coefficients in
$\GF(q)$, there exists a vector $u\in\ker (\widehat A(z^{-1}))$,
$u\neq 0$, with coordinates in the field $\GF(q^{r(z)})$. Thus
$f=z^\lambda\cdot u$ is a nonzero elementary solution with values
in the field $GF(q^{r(z)})$. Let $j\in\{1,\ldots,n\}$ be such that
$u_j\neq 0$.
 Then $f_j=u_j\cdot
z^\lambda\neq 0$. Moreover,
$$\Tr(f_j)=u_j\cdot z^\lambda+u_j^q\cdot (z^q)^\lambda +
\ldots+u_j^{q^{r(f)-1}}\cdot (z^{q^{r(f)-1}})^\lambda$$ is a
non-trivial linear combination of pairwise distinct characters.
Therefore $\Tr(f)$ is a nonzero solution  with values in the field
$GF(q)$, as required. \eproof

\brem\label{plaus} By Theorem 4.4 in \cite{Za2}, in the scalar
case the subspace $\ker (\Delta)\cap \cF_p(\Lambda,\GF(q))$ is
spanned by suitable traces of characters (i.e., of elementary
solutions). \erem

\section{Examples and applications}

The constructions in \ref{kuco}-\ref{grco} below are borrowed from
\cite{Ku1}.

\bsit\label{kuco} Let $\Lambda=\Z^s$ be a free abelian group
acting on a countable set $V$ with a finite quotient $V/\Lambda$
of cardinality $n$. Choosing a fundamental domain
$\Omega=\{v_1,\ldots,v_n\}\subseteq V$ (i.e., a set of
representatives of the orbits of $\Lambda$) we obtain a disjoint
partition
$$V=\bigsqcup_{v\in\Lambda} \tau_v(\Omega)\,.$$
According to this partition, every function $f\in \cF(V,\ka)$
yields a vector function $\vecti(f)\in\cF(\Lambda,\A^n_\ka)$. This
 provides an isomorphism of `fragmentation'
$$\vecti:\cF(V,\ka)\stackrel{\cong}\longrightarrow
\cF(\Lambda,\A^n)\,.$$ Any endomorphism $\Delta\in\End\,
(\cF(V,\ka))$ commuting with the $\Lambda$-action on the space
$\cF(V,\ka)$ amounts to an endomorphism $\Vecti(\Delta)\in\End\,
(\cF(\Lambda,\A^n))$ commuting with shifts. This yields a
representation
$$\Vecti:\End_\Lambda(\cF(V,\ka))\to\End_\Lambda\, (\cF(\Lambda,\A^n))\,,$$
where $\End_\Lambda$ stands for the set of endomorphisms commuting
with the action of $\Lambda$. In \ref{grco} and \ref{ovco} below
we provide two particular occasions of this construction.\esit

\bsit\label{grco} {\bf Covering graphs.} Let $\Gamma$ be a finite
connected graph with $n$ vertices. Consider the maximal abelian
covering $\widetilde\Gamma\to\Gamma$, which corresponds to the
commutator subgroup of the fundamental group $\pi_1(\G)$. The deck
transformation group (i.e., the Galois group) of the covering
$\widetilde\Gamma\to\Gamma$ is a free abelian group
$H_1(\Gamma;\Z)\cong\Z^{s}$, where $s=b_1(\Gamma)$. It acts freely
on $\tilde\G$ with quotient $\Gamma=\widetilde\Gamma/\Z^s$. Let
$V$ be the set of vertices of the periodic graph
$\widetilde\Gamma$. Every convolution operator $\Delta\in\End\,
(\cF(V,\ka))$  gives rise to an operator $\Vecti
(\Delta)\in\End_\Lambda(\cF(\Z^s,\A^n))$ commuting with shifts. In
particular, this applies  to the Markov and the Laplace operators
on $\widetilde\Gamma$ properly adopted to the case of positive
characteristic. As before, this yields a representation of the
algebra of convolution operators on $\tilde\G$. The same
construction works for any abelian covering of $\G$. \esit

To give further examples, let us consider the following toy model.

\bsit\label{ovco} {\bf Fragmentation: a toy model.} Let $\Lambda$
be a sublattice of finite index $n$ in a bigger lattice
$\widetilde\Lambda$. Consider the action of $\Lambda$ on the set
$V=\widetilde\Lambda$ by shifts. Then any endomorphism
$\Delta\in\End_\Lambda (\cF(\widetilde\Lambda,\ka))$ commuting
with shifts by $\Lambda$ gives rise to an endomorphism
$\Vecti_\Lambda(\Delta)\in\End_\Lambda (\cF(\Lambda,\A^n))$
commuting with shifts. In particular, for different choices of
$\Lambda$ this yields a collection of representations
$$\Vecti_\Lambda:\ka[\widetilde\Lambda]\to \End_\Lambda (\cF(\Lambda,\A^n))$$
of the convolution algebra on $\widetilde\Lambda$.\esit

New representations appear when restricting the morphism of
fragmentation to a translation invariant subspace e.g., to the
subspace $\cF_p(\Lambda,\A^n)\subseteq\cF(\Lambda,\A^n)$. We need
the following simple lemma.

\blem\label{repe}  In the notation as before, for any  $n$ coprime
to $p$ we have \be\label{vecti}
\vecti_\Lambda:\cF_p(\widetilde\Lambda,\ka)\stackrel{\cong}{\longrightarrow}
\cF_p(\Lambda,\A^n)\,.\ee \elem

\bproof For every $f=(f_1,\ldots,f_n)\in\cF_p(\Lambda,\ka)^n$, the
finite index sublattice  of $\Lambda$, $$\bigcap_{i=1}^n
\Lambda(f_i)\subseteq \Lambda(f)$$ is $p$-saturated. Therefore by
\ref{satu} $\Lambda(f)$ is $p$-saturated as well. This gives the
inclusion $\cF_p(\Lambda,\ka)^n\subseteq \cF_p(\Lambda,\A^n)$. The
opposite inclusion is evident. Now the last equality in
(\ref{vecti}) follows.

Letting furthermore $f=\vecti_\Lambda(\tilde f)$, where $\tilde
f\in\cF_p(\widetilde\Lambda,\ka)$, it is easily seen that the
sublattice of finite index
$$\Lambda(f)\supseteq \Lambda\cap\widetilde\Lambda (\tilde f)$$
is $p$-saturated. Thus $\vecti_\Lambda$ sends
$\cF_p(\widetilde\Lambda,\ka)$ injectively to
$\cF_p(\Lambda,\ka)^n$. In fact, this is a bijection. Indeed, let
$(v_1,\ldots,v_n)$ be a set of representatives of the cosets of
$\Lambda$ in $\widetilde\Lambda$.  For a vector function
$f=(f_1,\ldots,f_n)\in\cF_p(\Lambda,\ka)^n$ we let
$$\tilde f(v)=f_{i}(v-v_i)\quad\mbox{if}\quad v-v_i\in\Lambda\,.$$
It is easily seen that $f=\vecti_\Lambda(\tilde f)$. Moreover,
$\widetilde\Lambda (\tilde f)\supseteq \Lambda(f)$, where
$\Lambda(f)$ is $p$-saturated in $\widetilde\Lambda$. Hence
$\widetilde\Lambda (\tilde f)$ is also $p$-saturated and so
$\vecti_\Lambda$  in (\ref{vecti}) is an isomorphism.\eproof

In \ref{exo0}-\ref{exo3} below we apply our toy model \ref{ovco}
in different concrete settings.

\bexa\label{exo0} Consider a lattice $\widetilde \Lambda=\Z$ of
rank $s=1$ and the shift $\tau=\tau_{-1}$ on $\widetilde \Lambda$.
By Proposition \ref{sp1} the induced action of $\tau$ on the space
$\cF_p(\widetilde \Lambda,\ka)=\Span
(z^\lambda\,|\,z\in\ka^\times)$ can be diagonalized in the basis
of characters. More precisely, $\tau(z^\lambda)=z^{-1}\cdot
z^\lambda$, so $\spec (\tau)=\ka^\times$ and
$$\cF_p(\widetilde \Lambda,\ka)=\bigoplus_{\mu\in\ka^\times}
E_{\tau,\mu},\quad\mbox{where}\quad
E_{\tau,\mu}=\ka\cdot\mu^{-\lambda}\,.$$ The same is true for any
convolution operator $\Delta_a$ on $\widetilde \Lambda$, where
$a\in\ka[\widetilde \Lambda]$. Indeed, we have $\Delta_a=\widehat
a(\tau)$ and so
$$\Delta_a (z^\lambda)=\widehat
{a}(z^{-1})\cdot z^\lambda\quad\forall z\in\ka^\times\,.$$

Consider further a sublattice $\Lambda= n\Z$ of $\widetilde
\Lambda=\Z$ and the fragmentation
$\Vecti(\tau)=\Vecti_n(\tau)\in\End_\Lambda(\cF(\Lambda,\A^n))$.
Since $\tau(\delta_k)=\delta_{k+1}$ $\forall k\in\Z=\widetilde
\Lambda$, for the basic vector delta-functions
$\delta_{i,\lambda}=\delta_\lambda\cdot e_i$ on $\Lambda$ we have
$$\Vecti (\tau)
(\delta_{0,k})=\delta_{0,k+1}\,\,\,\mbox{for}\,\,\,k=0,\ldots,n-2,
\,\,\,\mbox{and}\,\,\, \Vecti (\tau)(\delta_{0,n-1})=
\delta_{1,0}\,.$$ The corresponding matrix of Laurent polynomials
$\widehat A(z)=\widehat {A_{\Vecti(\tau)}}(z)$ and its inverse
matrix $\widehat {A}(z)^{-1} =\widehat {A_{\Vecti(\tau^{-1})}}(z)$
are
$$\widehat A(z)=\left(%
\begin{array}{ccccc}
  0 & 0 & \ldots & 0 & z \\
  1 & 0 & \ldots & 0 & 0 \\
  0 & 1 & \ldots & 0 & 0 \\
  \vdots & \ddots & \ddots & \ddots & \vdots \\
  0 & 0 & \ldots & 1 & 0 \\
\end{array}%
\right)\,\,\,\mbox{resp.}\,\,\, \widehat {A}(z)^{-1}
=\left(%
\begin{array}{ccccc}
  0 & 1 & 0&  \ldots & 0  \\
  0 & 0 & 1 & \ldots & 0  \\
\vdots & \ddots & \ddots & \ddots & \vdots \\
  0 & 0 & \ldots & 0 & 1 \\

  z^{-1} & 0 & \ldots & 0 & 0 \\
\end{array}%
\right)=\widehat {A}^t(z^{-1})\,.$$ Thus
$$\det (xI_n-\widehat A(z))=x^n-z=\prod_{i=1}^n (x-\xi_i)\,,$$ where
$$\{\xi_1,\ldots,\xi_n\}=\{\sqrt[n]{z}\}=\spec(\widehat
{A}(z))\subseteq\ka^\times\,.$$  Suppose that $p\nmid n$. Then for
$z\in\ka^\times$ we have $\xi_i\neq\xi_j$ if $i\neq j$. So the
matrix $\widehat {A}(z^{-1})$ with eigenvalues
$\xi_1^{-1},\ldots,\xi_n^{-1}$ can be diagonalized in the basis of
eigenvectors
\be\label{eigen} v_i=\left(%
\begin{array}{c}
   1\\
  \xi_i \\
  \vdots \\
  \xi_i^{n-1} \\
\end{array}%
\right),\qquad i=1,\ldots,n\,,\ee where \be\label{eigen2}\widehat
{A}(z^{-1})(v_i)=\xi_i^{-1}v_i\quad\forall i=1,\ldots,n \,. \ee
The collection of elementary vector functions \be\label{eigen3}
\left(\lambda\longmapsto v_i\cdot z^\lambda\,|\,
 z\in\ka^\times,\,\,i=1,\ldots,n\right)\ee form a diagonalizing
basis for the endomorphism
$\Vecti_n(\tau)\in\End_\Lambda(\cF_p(\Lambda,\A^n))$ (cf.\
Proposition \ref{spde}(c)). Actually $v_i\cdot
z^\lambda=\vecti_n(\xi_i^\lambda)$, where $\lambda\longmapsto
\xi_i^\lambda$ is a character of $\tilde\Lambda=\Z$.

In conclusion, the isomorphism of fragmentation $\vecti_n$ as in
(\ref{vecti}) sends the basis of characters in
$\F_p(\widetilde\Lambda,\ka)$ to a basis of elementary functions
(\ref{eigen3}), which is diagonalizing for the endomorphism
$\Vecti_n(\tau)\in\End_\Lambda(\cF_p(\Lambda,\A^n))$. \eexa

\bexa\label{exo1} In the same setting as in \ref{exo0}, let us
consider an endomorphism
$$\Delta_{n,\phi}=\phi(\Vecti_n(\tau))=\Vecti_n(\phi(\tau))\in
\End_\Lambda (\cF(\Lambda,\A^n))$$ commuting with shifts, where
$\phi\in\ka[t,t^{-1}]$ is a Laurent polynomial. The matrix of
Laurent polynomials that corresponds to $\Delta_{n,\phi}$ is
$$\widehat A_{\Delta_{n,\phi}}=\phi(\widehat
A)=\widehat{\phi(A_{\Vecti_n(\tau)})} =
\widehat{A_{\Vecti_n(\phi(\tau))}}\,,$$ where $\widehat A$ is as
in \ref{exo0}. Recall that the equation $\Delta_{n,\phi}(f)=0$ for
$f\in\cF_p(\Lambda,\A^n)$ is equivalent to
$$\phi(\widehat{A}(z^{-1}))(\widehat{f}(z))=0\qquad\forall z\in
\ka^\times\,.$$
By the Spectral Mapping Theorem, for every
$z\in\ka^\times$ we have
$$\spec (\phi(\widehat {A}(z^{-1})))=\phi(\spec (\widehat
{A}(z^{-1})))=\{\phi(\xi_1^{-1}),\ldots,\phi(\xi_n^{-1})\}\,,$$
where as before $\{\xi_1,\ldots,\xi_n\}=\{\sqrt[n]{z}\}$. The
basis (\ref{eigen3}) is diagonalizing for $\Delta_{n,\phi}$.
Namely,
$$\Delta_{n,\phi} (v_i\cdot z^\lambda)=\phi(\xi_i^{-1})\cdot
v_i\cdot z^\lambda\,.$$ Indeed,
$$\widehat {A_{\phi(\Vecti_n(\tau))}}(z^{-1})\cdot\widehat {v_i\cdot
z^\lambda}= \widehat {A_{\phi(\Vecti_n(\tau))}}(z^{-1})\cdot
v_i\cdot\delta_z= \phi(\widehat
{A}(z^{-1}))(v_i)\cdot\delta_z$$$$=\phi(\xi_i^{-1})\cdot
v_i\cdot\delta_z=\phi(\xi_i^{-1})\cdot\widehat{v_i\cdot
z^\lambda}\,.$$ The Floquet-Fermi levels of $\Delta_{n,\phi}$ in
$\T^1=\ka^\times$ are
$$\Sigma_{\Delta_{n,\phi},\mu}=\{\xi^n\in\ka^\times\,|\,\phi(\xi^{-1})=\mu\}
=\iota\left(\phi^{-1}(\mu)\right)^n=\left(\phi^{-1}(\mu)\right)^{-n}\,.$$
In particular, the symbolic variety of $\Delta_{n,\phi}$ is
$$\Sigma_{\Delta_{n,\phi}}=\{\xi^n\in\ka^\times\,|\,\phi(\xi^{-1})=0\}\,.$$
Letting $\phi=q/z^\beta$, where $q=\prod_{j=1}^d
(t-\alpha_j)^{k_j}\in\ka[t]$ with $q(0)\neq 0$, we obtain
$$0\in\spec (\phi(\widehat
{A})(z^{-1}))\,\,\,\Longleftrightarrow\,\,\exists
i:\,\phi(\xi_i^{-1})=0 \,\,\,\Longleftrightarrow\,\,\exists
i:\,q(\xi_i^{-1})=0$$
$$\Longleftrightarrow\,\,\exists
i:\,\,\xi_i=\alpha_{j}^{-1}\,\,\,\mbox{for some}\,\,\,
j\in\{1,\ldots,d\}\,\,\,\Longleftrightarrow\,\,\,
z\in\{\alpha_1^{-n},\ldots,\alpha_d^{-n}\}\, ,$$ where as before
$\xi_i^n=z\in\ka^\times$, $i=1,\ldots,n$. Note that the sequence
$\alpha_j,\alpha_j^2,\ldots$ is periodic with period
$d_j=\ord_{\ka^\times}(\alpha_j)$ coprime to $p$.

The kernel $\ker (\Delta_{n,\phi})$ is spanned by elementary
vector functions $v_i\cdot z^\lambda$ as in (\ref{eigen3}) with
$z\in\{\alpha_1^{-n},\ldots,\alpha_d^{-n}\}$. Choose
$l\in\{0,\ldots,d_j-1\}$ such that $l\equiv n\mod d_i$. Letting
$z_{jl}=\alpha_{j}^{-l}=\xi_i^l$, such an elementary vector
function is proportional to the vector function
$$f_{jl}= z_{jl}^\lambda\cdot v_i
=(\quad\ldots,\, z_{jl}^{-1}v_i,\,v_i,\,
z_{jl}v_i,\,z_{jl}^2v_i,\,\ldots\quad)\in\cF_p(\Lambda,\A^n)\,,$$
where $v_i\in\A^n$ is an eigenvector of the matrix
$\widehat{A}(z_{jl}^{-1})$ as in (\ref{eigen}) with eigenvalue
$\alpha_{j}=\xi_i^{-1}$. Actually
$f_{jl}=\vecti_n(\alpha_{j}^{-\lambda})\,,$ where the character
$\alpha_{j}^{-\lambda}=\xi_i^\lambda:\tilde{\Lambda}\to\ka^\times$
of order $d_j$ (equal to its period) does not depend on $n$.\eexa

Let us illustrate the latter issue on a concrete example.

\bexa\label{exo2} Suppose that $p\neq 3$. Consider the equation
\be\label{crou}
\Delta_{n,\phi}(f)=\phi(\Vecti_n(\tau))(f)=\Vecti_n(\phi(\tau))(f)=0\,,\ee
where $\phi(t)=t+1+t^{-1}$. The roots $\alpha_1=\omega$ and
$\alpha_2=\omega^{-1}$ of $\phi$, where $\omega\in\ka^\times$ is a
primitive cubic root of unity, have multiplicative orders
$d_1=d_2=3$. There  exists an elementary solution $z^\lambda\cdot
v$ of (\ref{crou}) if and only if $z\in \{\omega^n,\omega^{-n}\}$,
respectively $\xi\in\{\omega^{-1},\omega\}$. Assuming that
$n\equiv 0\mod 3$ we obtain $l=0$ and $z=1$. Letting
$\xi=\omega^{-1}$ yields an elementary solution
$$f=(\ldots,v,v,v,\ldots),\quad\mbox{where}\quad v=\left(%
\begin{array}{c}
  1\\
  \omega^{-1} \\
  \omega \\
  1 \\
  \vdots\\
  \omega^{-1} \\
  \omega \\
\end{array}%
\right)\in\A^n\,.$$ This amounts to a $3$-periodic solution
$$\widetilde f=(\quad\ldots,\,1,\,\omega^{-1},\,\omega,\,1,\,\omega^{-1},
\,\omega,\,\ldots\quad)\in\cF_p(\tilde\Lambda,\ka)\,$$ of the
scalar equation \be\label{neweq}
\Delta_a(f)=\phi(\tau)(f)=0\quad\mbox{with}\quad
a=\delta_{-1}+\delta_0+\delta_{1}\in\ka[\tilde\Lambda]\,,\ee where
 as before $\tilde\Lambda=\Z$.
 The latter solution does not depend on $n$. Any
other elementary solution can be obtained from this one by a shift
and inversion $\iota:\omega\longmapsto\omega^{-1}$. Actually every
elementary solution of (\ref{neweq}) is proportional either to
$\widetilde f(\lambda)$ or to $\widetilde f(\lambda^{-1})$. Thus
the subspace of all $\Delta_a$-harmonic functions $$\ker
(\Delta_a)=\Span (\widetilde f,\widetilde
f\circ\iota)\subseteq\cF_p(\tilde\Lambda,\ka)$$ is
two-dimensional. It is easy to check that the same holds for every
$p\neq 3$ and $n\not\equiv 0\mod p$. \eexa

\bexa\label{exo3} For any $p\neq 2$ we consider the lattice
$\widetilde\Lambda=\Z^2$, its sublattice $\Lambda=2\Z\times
2\Z\subseteq \widetilde\Lambda$ of index 4,  and  the commuting
basic shifts $\tau_1=\tau_{(-1,0)}$ and $\tau_2=\tau_{(0,-1)}$ on
$\widetilde\Lambda$. In the standard basis $(e_{00},
e_{10},e_{01},e_{11})$ in $\A^4$ the corresponding matrices of
Laurent polynomials $\widehat{A_i}=\widehat{A_{\Vecti(\tau_i)}}$,
$i=1,2$, are as follows:
 $$\widehat{A_1}(z^{-1})
 =\left(%
\begin{array}{cccc}
  0 & z_1^{-1} & 0 & 0 \\
  1 & 0 & 0 & 0 \\
  0 & 0 & 0 & z_1^{-1} \\
  0 & 0 & 1 & 0 \\
\end{array}%
\right)\quad\mbox{resp.}\quad \widehat{A_2}(z^{-1})=
\left(%
\begin{array}{cccc}
  0 & 0 & z_2^{-1} & 0 \\
  0 & 0 & 0 & z_2^{-1} \\
  1 & 0 & 0 & 0 \\
  0 & 1 & 0 & 0 \\
\end{array}%
\right)\,.$$ In the common basis of eigenvectors
$$
v_{1}=\left(%
\begin{array}{c}
  1 \\
  x \\
  y \\
  xy \\
\end{array}%
\right),\quad
v_{2}=\left(%
\begin{array}{c}
  -1 \\
  x \\
  -y \\
  xy \\
\end{array}%
\right),\quad
v_{3}=\left(%
\begin{array}{c}
  -1 \\
  -x \\
  y \\
  xy \\
\end{array}%
\right),\quad
v_{4}=\left(%
\begin{array}{c}
  1 \\
  -x \\
  -y \\
  xy \\
\end{array}%
\right)\,$$ these commuting matrices are reduced to the diagonal
form
$$\diag \left(x^{-1},\,-x^{-1},\, x^{-1},\,-x^{-1}\right) \quad\mbox{resp.}\quad
\diag\left(y^{-1},\,y^{-1},\,-y^{-1},\,-y^{-1}\right)\,,$$ where
$x^2=z_1$ and $y^2=z_2$.

 Consider further an endomorphism $\rho$ commuting
with shifts, where
$$\rho=\phi(\Vecti(\tau_1),\Vecti(\tau_2))=
\Vecti(\phi(\tau_1,\tau_2))\in\End_\Lambda (\cF(\Lambda,\A^4))\,$$
and $\phi\in\ka[z_1,z_2, z_1^{-1},z_2^{-1}]$ is a Laurent
polynomial. The same argument as in \ref{exo1} above shows that
the basis of elementary vector functions
$$\left(\lambda\longmapsto v_j(z)\cdot z^\lambda\,
|\,z\in \T^2,\,\,j=1,2,3,4\right)\,$$ is diagonalizing for $\rho$.
More precisely, letting
$$\phi_1(x,y)=\phi(x,y),\quad \phi_2(x,y)=\phi(-x,y),
\quad\phi_3(x,y)=\phi(x,-y),\quad\phi_4(x,y)=\phi(-x,-y)\,,$$ by
the Spectral Mapping Theorem we obtain
$$\rho(v_j\cdot z^\lambda)=\phi_j(x^{-1},y^{-1})\cdot v_j\cdot z^\lambda
\,
$$
(cf.\ Theorem \ref{spde}(c)). The  Floquet-Fermi level-$\mu$ set
of $\rho$ is
$$\Sigma_{\rho,\mu}=\{z=(x^2,y^2)\in\T^2\,|\,\exists
j\in\{1,2,3,4\}: \phi_j(x^{-1}, y^{-1})=\mu\}\,$$ (see \ref{jd}).
The eigenspace of $\rho$ with eigenvalue $\mu$ is
$$E_{\rho,\mu}=\Span \left(v_j\cdot z^\lambda\,
|\,z=(x^2,y^2)\in \T^2: \phi_j(x^{-1}, y^{-1})=\mu,\,
j=1,2,3,4\right)\,.$$ Moreover,
$$\cF_p(\Lambda,\A^4)=\bigoplus_{\mu\in\ka} E_{\rho,\mu}\,.$$
In particular,
the symbolic hypersurface of $\rho$ is
$$\Sigma_{\rho}=\{z=(x^2,y^2)\in\T^2\,|\, \exists
j\in\{1,2,3,4\}:\phi_j(x^{-1}, y^{-1})=0\}\,.$$ The space
$E_{\rho,0}$ of all $\rho$-harmonic vector functions is
$$\ker (\rho)=\Span \left(v_j\cdot z^\lambda\,
|\,z=(x^2,y^2)\in \T^2:\phi_j(x^{-1},
y^{-1})=0,\,\,j=1,2,3,4\right)\,.$$ \eexa

\end{document}